\documentclass[12pt]{article}
\usepackage{mathtools}
\usepackage[utf8]{inputenc}
\usepackage{amssymb,amsmath,amsfonts,eucal,mathrsfs,amsthm}
\usepackage[colorinlistoftodos]{todonotes}
\usepackage{xcolor}
\usepackage{enumitem}
\newtheorem{theorem}{Theorem}
\newtheorem{proposition}[theorem]{Proposition}
\newtheorem{lemma}[theorem]{Lemma}

\newtheorem{example}[theorem]{Example}

\theoremstyle{definition}
\newcommand{\R}{\mathbb{R}}
\newcommand{\Z}{\mathbb{Z}}
\newcommand{\Q}{\mathbb{Q}}
\newcommand{\Sf}{\mathbb{S}}

\newcommand{\spa}{\mbox{span}}

\newcommand{\Ric}{\mbox{Ric}}

\def\<{{\langle}}
\def\>{{\rangle}}
\def\B{\mathcal{B}}
\def\CP{\mathord{\mathbb C}\mathord{\mathbb P}}

\def\n{\nabla}

\def\a{\alpha}

\def\be{\begin{equation} }
\def\ee{\end{equation} }

\begin{document}

\title{Geometric and topological rigidity of pinched submanifolds II}
\author{Theodoros Vlachos}
\date{}
\maketitle

\renewcommand{\thefootnote}{\fnsymbol{footnote}} 
\footnotetext{\emph{2020 Mathematics Subject Classification.} 53C40, 53C42.} 
\renewcommand{\thefootnote}{\arabic{footnote}} 

\renewcommand{\thefootnote}{\fnsymbol{footnote}} 
\footnotetext{\emph{Keywords and phrases.} Homology groups, pinching, 
length of the second fundamental form, sectional, isotropic and 
mean curvature, Bochner-Weitzenb\"ock operator.} 
\renewcommand{\thefootnote}{\arabic{footnote}}

\begin{abstract}
We continue the study of the geometry and topology of
compact submanifolds of arbitrary codimension in space
forms that satisfy a pinching condition involving the length
of the second fundamental form and the mean curvature.
Our primary focus is on four-dimensional submanifolds,
where both the results obtained and the methods employed
differ substantially and are considerably more intricate than
in higher dimensions. This study relies critically on concepts
from four-dimensional geometry, the theory of Riemannian
manifolds with nonnegative isotropic curvature, and the
Bochner technique, each playing an essential role. The results
are sharp and extend previous work by several authors,
without imposing additional assumptions on either the mean
curvature or the fundamental group of the submanifold.
\end{abstract}

\section{Introduction}
A fundamental problem in differential geometry is 
to understand the interplay between the geometry 
and topology of Riemannian manifolds. In the context 
of submanifold theory, it is particularly 
interesting to investigate how the geometry and topology 
of submanifolds in space forms are influenced by 
pinching conditions involving either intrinsic or 
extrinsic curvature invariants.

For minimal submanifolds of spheres with a sufficiently
pinched second fundamental form, a fundamental result was
first established by Simons in his seminal work
\cite{Sim}. Later, Chern, do Carmo, and Kobayashi
\cite{CdCK} proved a celebrated rigidity theorem.
Their contributions have since inspired a wealth of
important developments in the study of pinching conditions,
including, for example, \cite{adc, HW, Le1, OV, SX, V2, XG}.

In our previous paper \cite{v3}, we studied the
geometric and topological rigidity of submanifolds
$f\colon M^n \to \Q_c^{n+m}$ satisfying the pointwise
pinching condition
\be\label{1}
S \leq a(n,k,H,c) \tag{$\ast$}
\ee
where $S$ denotes the \emph{squared length} of the second
fundamental form
$\alpha_f \colon TM \times TM \to N_fM$,
taking values in the normal bundle $N_fM$.
The \emph{mean curvature} is defined as
$H = \| \mathcal H \|$, the norm of the
\emph{mean curvature vector field}
$$
\mathcal H = \frac{1}{n}\,\mathrm{tr}\, \alpha_f,
$$ 
where $\mathrm{tr}$ denotes the trace.
Here $k$ is an integer with $1 \leq k \leq n/2$,
and the function $a$ is defined as
$$
a(n,k,t,c)
=nc+\frac{n^3t^2}{2k(n-k)}-\frac{n|n-2k|}{2k(n-k)}t\sqrt{
n^2t^2+4ck(n-k)},\;t, c\geq0.
$$
The ambient space $\Q_c^{n+m}$ is the 
$(n+m)$-dimensional, complete, simply connected 
space form of constant curvature $c$. For simplicity, 
we restrict to the cases $c\in\{0,1\}$, unless 
stated otherwise. Accordingly, $\Q_c^{n+m}$ is 
either the Euclidean space $\R^{n+m}$ when $c=0$, 
or the unit sphere $\Sf^{n+m}$ when $c=1$.

A straightforward computation shows that the standard 
embedding of a torus 
$$
\mathbb T^n_k(r)=\mathbb S^k(r)\times
\mathbb S^{n-k}(\sqrt{1-r^2}),\;1\leq k\leq n-1,
$$
into the unit sphere $\mathbb S^{n+1}$, where 
$\mathbb S^k(r)$ denotes the $k$-dimensional 
sphere of radius $r<1$, satisfies 
\eqref{1} with equality if $r\geq\sqrt{k/n}$, 
or $n=2k$. In all other cases, one has $S>a(n,k,H,1)$.

The pinching condition \eqref{1} has been studied
primarily in the case of the smallest admissible value
of $k$ (see, e.g., \cite{WXia,Xu,XG,ZZ}). 
Shiohama and Xu \cite{SX} proved that compact
submanifolds in space forms of nonnegative curvature
are homeomorphic to a sphere, provided \eqref{1} holds
as a strict inequality at every point when $k=1$. In \cite{V2}, a 
homology vanishing result was established for 
submanifolds in space forms of nonnegative 
curvature that satisfy the strict pointwise inequality in 
\eqref{1} for some integer $1\leq k\leq n-1$. 
The argument there relied on the nonexistence of stable
currents, as established by Lawson and Simons \cite{LS},
under suitable upper bounds for the second
fundamental form.

In our recent work \cite{v3}, we studied the geometric
and topological rigidity of submanifolds of dimension
$n \geq 5$ satisfying the pinching condition \eqref{1},
without assuming strictness or imposing any restrictions
on the mean curvature.
Among other results, we established the following theorem, 
showing that the
pinching condition either forces the vanishing of
homology in certain intermediate dimensions or uniquely
determines the submanifold up to congruence.
Recall that a submanifold is said to be substantial if
it is not contained in any proper totally geodesic
submanifold of the ambient space.

\begin{theorem}\label{th1}
Let $f\colon M^n\to\Q_c^{n+m}$, 
$n\geq 5,c\geq0$, be a substantial isometric 
immersion of a compact oriented Riemannian 
manifold. Assume that inequality \eqref{1} is 
satisfied for some integer integer $2\leq k\leq n/2$ 
at every point. Then, either the homology groups of 
$M^n$ satisfy 
$$
H_p(M^n;\Z)=0\;\,\text{for all}\;\,k\leq p\leq n-k
\;\,\text{and}\;\, H_{k-1}(M^n;\Z)\cong\Z^{\beta_{k-1}(M)},
$$
where $\beta_{k-1}(M)$ 
is the $(k-1)$-th Betti number, or equality holds 
in \eqref{1} at every point and one of the following 
assertions holds:
\vspace{1ex}

\noindent $(i)$ $M^n$ is isometric to a Clifford 
torus $\mathbb T^n_p(\sqrt{p/n}),k\leq p\leq n/2$, 
and $f$ is the standard minimal embedding into 
$\Sf^{n+1}$.
\vspace{1ex}

\noindent $(ii)$ $M^n$ is isometric 
to a torus $\mathbb T^n_k(r)$ with $r>\sqrt{k/n}$ 
and $f$ is the standard embedding into $\Sf^{n+1}$.
\vspace{1ex}

\noindent $(iii)$ $M^n$ is isometric 
to a torus $\Sf^k(r)\times \Sf^k(\sqrt{R^2-r^2})$ 
and $f$ is a composition $f=j\circ g$, where 
$g\colon M^n\to\Sf^{n+1}(R)$ is the standard 
embedding of the torus 
$\Sf^k(r)\times \Sf^k(\sqrt{R^2-r^2})$ into a 
sphere $\Sf^{n+1}(R)$ and 
$j\colon\Sf^{n+1}(R)\to\Q_c^{n+2}$ is an 
umbilical inclusion with $c=0,1$ and $R<1$ if $c=1$.
\end{theorem}

The main focus of this paper is the study of four-dimensional
submanifolds satisfying condition \eqref{1}, which are
not covered by Theorem \ref{th1}.
The results obtained here, as well as the methods
employed, differ substantially from those in the proof
of Theorem \ref{th1} for higher dimensions and are
considerably more intricate.
Our approach relies essentially on tools from
four-dimensional geometry, the theory of Riemannian
manifolds with nonnegative isotropic curvature, and the
Bochner technique.
A key observation is that four-dimensional submanifolds
satisfying the pinching condition \eqref{1} necessarily
have nonnegative isotropic curvature, a concept
introduced and studied by Micallef and Moore
\cite{MM}.

Throughout the paper, all submanifolds under 
consideration are assumed to be connected. 
The main result of the paper can be stated 
as follows.

\begin{theorem}\label{k=2}
Let $f\colon M^4\to\Q_c^{4+m},c\geq0$, be a 
substantial isometric immersion of a compact, 
oriented Riemannian four-manifold. Suppose 
inequality \eqref{1} is satisfied for $k=2$ 
at every point. Then, one of the following 
assertions holds:
\vspace{1ex}

\noindent $(i)$ $M^4$ is diffeomorphic 
to $\Sf^4$.  
\vspace{1ex}

\noindent $(ii)$ The universal cover of $M^4$ is 
isometric to a Riemannian product 
$\R\times N$, where $N$ is diffeomorphic to 
$\Sf^3$ with nonnegative Ricci curvature.  
\vspace{1ex}

\noindent $(iii)$
Equality holds in \eqref{1} for $k=2$ at every point, 
and one of the following holds: 
\vspace{0.5ex}
\begin{enumerate}[topsep=1pt,itemsep=1pt,partopsep=1ex,parsep=0.5ex,leftmargin=*, label=(\roman*), align=left, labelsep=-0.5em]
\item [(a)] $M^4$ is isometric to a torus 
$\Sf^2(r)\times \Sf^2(\sqrt{R^2-r^2})$ and 
$f$ is a composition $f=j\circ g$, where 
$g\colon M^4\to\Sf^5(R)$ is the standard 
embedding of the torus 
$\Sf^2(r)\times \Sf^2(\sqrt{R^2-r^2})$ into a 
sphere $\Sf^5(R)$, and $j\colon\Sf^5(R)\to\Q_c^6$ 
is an umbilical inclusion with $R\leq1$ if $c=1$.
\item [(b)] $M^4$ is isometric to the 
projective plane $\CP_r^2$ of constant holomorphic 
curvature $4/3r^2$ with $r=1/\sqrt{c+H^2}$ and 
$f=j\circ g$, where $g$  
is the standard minimal embedding 
of $\CP_r^2$ into $\Sf^7(r)$, and 
$j\colon\Sf^7(r)\to\Q_c^8$ is an umbilical inclusion.
\end{enumerate}
\end{theorem}

In the final section of the paper, we present a 
method for constructing geometrically distinct 
isometric immersions of manifolds diffeomorphic 
either to $\Sf^4$ or to the torus 
$\Sf^1\times\Sf^3$, which also satisfy 
\eqref{1} for $k=2$. These examples 
demonstrate that Theorem \ref{k=2} is 
sharp. In particular, the standard 
isometric embedding of the torus 
$\Sf^1(r)\times\Sf^3(\sqrt{1-r^2})$ into $\Sf^5$ 
satisfies \eqref{1} for $k=2$ whenever  
$r\geq1/2$, and its fundamental 
group is infinite. This example shows 
that the class of submanifolds appearing 
in part $(ii)$ of Theorem \ref{k=2} is nonempty. 

For submanifolds of arbitrary dimension 
satisfying the pinching condition \eqref{1} 
with $k=2$, we establish the following result.

\begin{theorem}\label{ck=2}
Let $f\colon M^n\to\Q_c^{n+m}$, with
$n\geq 4$ and $c\geq0$, be a substantial isometric 
immersion of a compact, oriented Riemannian 
manifold. Assume that inequality \eqref{1} is 
satisfied at every point for $k=2$. 
Then, one of the following holds: 
\vspace{1ex}

\noindent $(i)$ The manifold $M^n$ is diffeomorphic to $\Sf^4$.
\vspace{1ex}

\noindent $(ii)$ The universal cover 
of $M^n$ is isometric to a Riemannian product 
$\R\times N$, where $N$ is diffeomorphic to 
$\Sf^3$ with nonnegative Ricci curvature. 
\vspace{1ex}

\noindent $(iii)$ If $n\geq5$, the homology of $M^n$ 
satisfies 
$$
H_p(M^n;\Z)=0\;\;\text{for all}\;\;2\leq p\leq n-2
$$
and
$$
H_p(M^n;\Z)\cong\Z^{\beta_1(M)}\;\;\text{for}\;\;p=1,n-1.
$$
\noindent $(iv)$ Equality holds in condition \eqref{1} for $k=2$ 
at every point, and the submanifold is as described in 
parts $(i)$-$(iii)$ of Theorem \ref{th1}, 
or as in part $(iii)$ of Theorem \ref{k=2}. 
\end{theorem}

The results presented here extend those of previous studies 
\cite{adc, AC, CdCK, HW, Le1, Sim, XG, XZ, ZZ} without 
imposing additional assumptions on either the mean 
curvature or the fundamental group of the submanifold. 
In the final section, we further construct geometrically 
distinct isometric immersions of manifolds diffeomorphic 
to either the sphere $\Sf^n$ or the torus $\Sf^1 \times\Sf^{n-1}$, 
with $n\geq 4$, which also satisfy \eqref{1} for $c=0$ and 
$k=2$. These examples illustrate the sharpness of 
Theorem \ref{ck=2}. Finally, note that the homology 
described in part $(iii)$ of Theorem \ref{ck=2} coincides 
with that of the sphere $\Sf^n$ with $g=\beta_1(M)$ 
handles attached, or equivalently, with the homology 
of the connected sum $\#_g(\Sf^1 \times\Sf^{n-1})$ 
consisting of $g$ summands.

The geometric and topological rigidity of submanifolds of 
arbitrary dimension satisfying the pinching 
condition \eqref{1} for $k=1$ was studied in \cite{v3}. Combined  
with the results presented here, this provides a comprehensive 
understanding of compact submanifolds satisfying \eqref{1} for the 
two smallest admissible values of the integer $k$.

It is worth noting that in the recent paper \cite{a}, the 
authors extended the pinching result of Simons \cite{Sim} 
and that of Chern, do Carmo, and Kobayashi 
\cite{CdCK} to complete (not necessarily compact) 
minimal submanifolds in spheres. It would be 
of interest to extend the results of the present 
paper, as well as those in \cite{v3}, to complete 
submanifolds in space forms, without imposing 
any restrictions on the mean curvature. Such a 
generalization could provide deeper insights into 
the geometry and topology of submanifolds satisfying 
the pinching condition \eqref{1}, beyond the minimal or 
compact setting.

\section{Background}
\subsection{Geometry of 4-dimensional manifolds}
In this section, we collect basic facts about four-dimensional 
geometry. For a detailed exposition of the subject, we refer 
the reader to \cite{LeB, sea}.

Let $(M,\<\cdot,\cdot\>)$ be an oriented Riemannian 
manifold of dimension $n=4$ with Levi-Civit\'a connection 
$\n$ and curvature tensor $R$ given by
$$
R(X,Y)=[\n_X,\n_Y]-\n_{[X,Y]},\;X,Y\in\mathcal X(M).
$$
The Ricci tensor of $(M,\<\cdot,\cdot\>)$ is defined by
$$
\Ric(X,Y)=\sum_i\<R(X,E_i)E_i,Y\>,\;X,Y\in\mathcal X(M),
$$
where $\{E_i\}_{1\leq i\leq 4}$ is a local orthonormal frame.

At any point $x\in M$, we consider the {\emph 
{Bochner-Weitzenb\"ock operator} $\B^{[2]}$ as an 
endomorphism of the space of 2-vectors 
$\Lambda^2T_xM$ at $x$ given by 
\begin{align}\label{WB}
\<\<\B^{[2]}(v_1\wedge v_2), w_1\wedge w_2\>\>
=&\;\Ric(v_1,w_1)\<v_2,w_2\>+\Ric(v_2,w_2)\<v_1,w_1\>\nonumber\\
&-\Ric(v_1,w_2)\<v_2,w_1\>-\Ric(v_2,w_1)\<v_1,w_2\>\nonumber\\
&-2\<R(v_1,v_2)w_2,w_1\>,
\end{align}
and then extend it linearly to all of $\Lambda^2T_xM$.
Here $\<\<\cdot,\cdot\>\>$ stands for the inner product of 
$\Lambda^2T_xM$ defined by 
$$
\<\<v_1\wedge v_2,w_1\wedge w_2\>\>=\det(\<v_i,w_j\>).
$$
The Bochner-Weitzenb\"ock operator is a self-adjoint operator. 
If $Z\in\Lambda^2T_xM$, the dual 2-form $\omega$ is 
defined by 
$$
\omega(v,w)=\<\<Z,v\wedge w\>\>, 
$$ 
and we may regard $Z$ as the skew-symmetric 
endomorphism of the tangent space at $x$ via  
$$
\<Z(v),w\>=\<\<Z,v\wedge w\>\>. 
$$

Clearly $\B^{[2]}$ can also be viewed as an endomorphism 
of the bundle $\Omega^2(M)$ of 2-forms of the manifold via 
the inner product $\<\<\cdot,\cdot\>\>$. If $\omega$ is a 
2-form, then $\B^{[2]}(\omega)$ is given by 
$$
\B^{[2]}(\omega)(X_1,X_2)=\omega(\Ric(X_1),X_2)
+\omega(X_1,\Ric(X_2))-\sum_i\omega(R(X_1,X_2)E_i,E_i).
$$
Then the Bochner-Weitzenb\"ock operator acts on 
$\Lambda^2TM$ by 
$$
\omega(\B^{[2]}(X_1\wedge X_2))
=\B^{[2]}(\omega)(X_1,X_2).
$$ 
Taking $\omega$ to 
be the dual form to $w_1\wedge w_2$ yields \eqref{WB}. 
Throughout the paper, we will generally identify 2-forms 
with 2-vectors.

We recall the Bochner-Weitzenb\"ock formula which 
can be written as 
$$
\<\Delta\omega,\omega\>=\frac{1}{2}\Delta\|\omega\|^2+\|\n\omega\|^2
+\<\B^{[2]}(\omega),\omega\>
$$
for any $\omega\in\Omega^2(M)$. This implies that any 
harmonic 2-form on a compact manifold is parallel provided 
that the Bochner-Weitzenb\"ock operator is nonnegative. 

The bundle of 2-forms of any oriented 4-dimensional 
Riemannian manifold $M$ decomposes as a direct sum 
$$
\Omega^2(M)=\Omega_+^2(M)\oplus\Omega_-^2(M)
$$
of the eigenspaces of the Hodge star operator 
$$
\ast\colon\Omega^2(M)\to\Omega^2(M)
$$ 
The sections of 
$\Omega_+^2(M)$ are called {\emph {self-dual 2-forms}}, 
whereas the ones of $\Omega_-^2(M)$ are called 
{\emph {anti-self-dual 2-forms}}. 
Accordingly, at any point $x\in M$, we have the splitting
$$
\Lambda^2T_xM=\Lambda_+^2T_xM\oplus\Lambda_-^2T_xM, 
$$
where $\Lambda_\pm^2T_xM$ 
are the eigenspaces of the Hodge star operator 
$$
\ast\colon\Lambda^2T_xM\to\Lambda^2T_xM.
$$ 
Both spaces $\Lambda_+^2T_xM$ and $\Lambda_-^2T_xM$ 
are $\B^{[2]}$-invariant (see Proposition 1 in \cite{smichig}). 
Hence, we obtain the corresponding decomposition 
$$
\B^{[2]}=\B_+^{[2]}\oplus\B_-^{[2]}. 
$$

Suppose now that $M$ is a compact, 
oriented Riemannian four-manifold. The Hodge theorem 
guarantees that every de Rham cohomology class on $M$ 
has a unique harmonic representative. In particular, the space 
$\mathscr H^2(M)$ of harmonic 2-forms decomposes as 
$$
\mathscr H^2(M)=\mathscr H^2_+(M)\oplus\mathscr H^2_-(M),
$$ 
where $\mathscr H^2_+(M)$ and $\mathscr H^2_-(M)$ are 
the spaces of self-dual and anti-self-dual harmonic 2-forms, 
respectively. The dimensions of these subspaces, denoted by
$\beta_\pm(M)=\dim\mathscr H^2_\pm(M)$, are oriented 
homotopy invariants of $M$. Their difference 
$$
\sigma=\beta_+(M)-\beta_-(M)
$$ 
is the signature of $M$, 
while their sum equals the second Betti number 
$\beta_2(M)$ of the manifold. 

\subsection{Isotropic curvature}
Let $(M,g)$ be a Riemannian manifold of dimension $n\geq4$. 
We say that $(M,g)$ has \emph {nonnegative isotropic curvature} 
at a point $x\in M$ if
$$
R_{1331} + R_{1441}+ R_{2332} + R_{2442}
-2 R_{1234} \geq0, 
$$
for all orthonormal four-frames 
$\{e_1,e_2,e_3,e_4\}\subset T_xM$. Here, $R_{ijk\ell}$ are the 
components of the curvature tensor $R$, defined by 
$$
R_{ijk\ell}=g\big(R(e_i,e_j)e_k, e_\ell\big), \;1\leq i,j,k,\ell\leq4. 
$$ 
If the strict inequality holds, we say that $(M,g)$ has 
\emph{positive isotropic curvature} at $x$. 
The manifold $(M,g)$ is said to have nonnegative (or positive) 
isotropic curvature if it satisfies the corresponding condition 
at every point and for all orthonormal four-frames.

The following result is well known (see, for instance, \cite{MWolf}).

\begin{lemma}\label{sb}
For any four-dimensional Riemannian manifold $M$, the 
non-negativity of the isotropic curvature at a point $x\in M$ 
is equivalent to the non-negativity of the Bochner-Weitzenb\"ock 
operator $\B^{[2]}$ at $x$. 
\end{lemma}

\subsection{The pinching condition and the Lawson-Simons inequality}
Lawson and Simons \cite{LS} proved that specific 
bounds on the second fundamental form for submanifolds 
of spheres can result to the vanishing of homology groups 
with integer coefficients. By employing the second variation 
of area, they effectively ruled out stable minimal currents in 
certain dimensions. Consequently, since the existences theorems 
of Federer and Fleming \cite{FF} allow area minimization 
within a homology class, this result implies 
the triviality of integral homology.

The result of Lawson and Simons \cite{LS} mentioned above 
was later strengthened by Elworthy and Rosenberg 
\cite[p.\ 71]{ER} without requiring the bound on the second 
fundamental form to be strict at all points of the submanifold. 
In this section, we state their theorem and then quote an 
auxiliary result from \cite{v3} that establishes the relation 
between our pinching condition \eqref{1} and the inequality 
\eqref{LS} below required in their result.

\begin{theorem}(\cite{ER, LS, X})\label{er} Let 
$f\colon M^n\to\Q_c^{n+m}, n\geq 4, c\geq0$, 
be an isometric immersion of a compact Riemannian 
manifold. Let $p$ be an integer such that $1\leq p\leq n-1$. 
Assume that at any point $x\in M^n$ and for every 
orthonormal basis $\{e_1,\ldots,e_n\}$ of $T_xM$, 
the second fundamental form 
$\alpha_f\colon TM\times TM\to N_fM$ satisfies
\be\label{LS}
\sum_{i=1}^p\sum_{j=p+1}^n\big(2\|\a_f(e_i,e_j)\|^2
-\<\a_f(e_i,e_i),\a_f(e_j,e_j)\>\big)\leq p(n-p)c \tag{$\ast\ast$}. 
\ee
If there exists a point at which inequality \eqref{LS} is 
strict for every orthonormal basis of the tangent space 
at that point, then the homology groups satisfy
$$
H_p(M^n;\mathbb{Z})=H_{n-p}(M^n;\mathbb{Z})=0.
$$
\end{theorem} 

Next, we recall the following result from \cite{v3}, which 
establishes the relationship between the inequalities 
\eqref{1} and \eqref{LS} for four-dimensional submanifolds.

\begin{lemma}\label{propu}
Let $f\colon M^4\to\Q_c^{4+m},c\geq0$, be an 
isometric immersion of a $4$-dimensional Riemannian 
manifold $M^4$ such that inequality \eqref{1} is satisfied at 
a point $x\in M^4$ for $k=2$. Then the following 
assertions hold at $x$:
\vspace{1ex}

\noindent $(i)$ The inequality \eqref{LS} holds for $p=2$ 
and for all orthonormal bases of $T_xM$. Moreover, 
if inequality \eqref{1} is strict at $x$, then \eqref{LS} 
is also strict for every orthonormal basis of $T_xM$.
\vspace{1ex}

\noindent $(ii)$ Suppose now that equality holds in \eqref{LS} 
for some orthonormal basis $\{e_i\}_{1\leq i\leq 4}$ 
of $T_xM$ and $p=2$. Then there exist normal vectors 
$\eta_1, \eta_2 \in N_fM(x)$ such that the shape 
operator $A_\xi$ associated to any $\xi\in N_fM(x)$ 
satisfies 
\be\label{Vj}
\pi_j\circ A_\xi|_{V_j}=\<\xi,\eta_j\>Id,\; j=1,2, \nonumber
\ee
where $Id$ is the identity map on the tangent space at $x$, 
$V_1=\spa\{e_1,e_2\}$, $V_2=\spa\{e_3,e_4\}$ and $\pi_j$ 
denotes the projection onto $V_j,j=1,2$. 
\end{lemma}

\section{Proofs of the main results}

We recall a key result, namely Proposition 16 from \cite{OV}, 
which provides an estimate for the Bochner-Weitzenb\"ock 
operator of a submanifold in terms of its second fundamental 
form. For four-dimensional submanifolds, this proposition 
can be stated as follows:

\begin{proposition}\label{lemp}
Let $f\colon M^4\to\Q_c^{4+m}$ be an isometric 
immersion of a $4$-dimensional manifold $M^4$. The 
Bochner-Weitzenb\"ock operator $\B^{[2]}$ of 
$M^4$ satisfies the following pointwise inequality 
\be\label{ineqp}
\mathop{\min_{\omega\in \Omega^2(M^4)}}_{\|\omega\|=1}\<\B^{[2]}\omega,\omega\>
\geq 4c+8H^2-S.
\ee
If equality holds in \eqref{ineqp} at a 
point $x\in M^4$, then for every unit vector $\xi \in N_f M(x)$,
the shape operator $A_\xi(x)$ has at most two distinct 
eigenvalues, each of multiplicity $2$. 
\end{proposition}

We now state the following auxiliary results.

\begin{proposition}\label{nnic}
Let $f\colon M^4\to\Q^{4+m}_c, c\geq 0$, be an 
isometric immersion of an oriented $4$-dimensional 
Riemannian manifold. Suppose that $f$ satisfies inequality 
\eqref{1} at every point for $k=2$. Then, the 
following assertions hold:\vspace{1ex}

\noindent $(i)$ The Bochner-Weitzenb\"ock operator 
$\B^{[2]}$ is nonnegative, and the manifold $M^4$ 
consequently has nonnegative isotropic curvature.
\vspace{1ex}

\noindent $(ii)$ Suppose that $M^4$ is compact and 
satisfies $\beta_2(M)\neq0$. Then equality holds in \eqref{1} 
at every point, and the shape operator, associated 
with any normal direction has at most two distinct 
eigenvalues, each of multiplicity $2$. Moreover, 
at each point $x\in M^4$, there exists an oriented 
orthonormal four-frame 
$\{e_1,e_2,e_3,e_4\}\subset T_xM$ such that 
the second fundamental form 
$\a_f$ of $f$ satisfies the following conditions: 
\begin{align}
\;\;\;\;\a_{11}=\a_{22},\,\a_{33}=&\a_{44},\,\a_{12}=\a_{34}=0,
\,\|\a_{23}\|=\|\a_{14}\|,\,\|\a_{24}\|=\|\a_{13}\|,\label{a1}\\
\<\a_{14}+\a_{23}&,\a_{13}-\a_{24}\>=0,\;\,\<\a_{13}+\a_{24},\a_{14}-\a_{23}\>=0,\label{a2}\\
&\|\a_{13}\|^2+\|\a_{14}\|^2=c+\<\a_{11},\a_{44}\>,\label{a3}
\end{align}
where, for simplicity, we denote $\a_{ij}=\a_f(e_i,e_j)$.
\end{proposition}

\proof \noindent $(i)$ From \eqref{1}, if follows 
that the right-hand side of inequality \eqref{ineqp} in 
Proposition \ref{lemp} is nonnegative. Consequently, the 
Bochner-Weitzenb\"ock operator $\B^{[2]}$ is 
nonnegative. By Lemma \ref{sb}, this implies that 
the manifold $M^4$ has nonnegative isotropic curvature. 
\vspace{0.5ex}

\indent $(ii)$ By the Hodge theorem, there exists a 
nontrivial harmonic 2-form $\omega\in\Omega^2(M)$. 
The Bochner-Weitzenb\"ock formula then implies 
$$
\<\B^{[2]}\omega,\omega\>=0
$$ 
at every point. Therefore, equality holds in inequality 
\eqref{ineqp}, and as a result, Proposition \ref{lemp} 
implies that the shape operator, associated with any 
normal direction has at most two distinct 
eigenvalues, each of multiplicity $2$. 

On the other hand, from part $(i)$ of Lemma \ref{propu}, 
we know that inequality \eqref{LS} holds for $p=2$ 
at every point and for every orthonormal basis. Since 
$H_2(M;\mathbb{Z})\neq0$, Theorem \ref{er} implies 
that at each point $x\in M^4$ there exists an oriented 
orthonormal four-frame $\{e_1,e_2,e_3,e_4\}\subset T_xM$ 
such that equality holds in \eqref{LS} 
for $p=2$. It follows from part $(i)$ of Lemma \ref{propu} 
that equality also holds in \eqref{1} at any point. Furthermore, 
by choosing an orthonormal normal basis
$\{\xi_\alpha\}_{1\leq\alpha\leq m}$ at $x\in M^4$ such 
that the mean curvature vector is $\mathcal{H}(x)=H\xi_1$, 
part $(ii)$ of Lemma \ref{propu} implies that the 
corresponding shape operators $A_{\xi_\a},1\leq\a\leq m$, 
are given by 
\be\label{A}
\begin{cases}
A_{\xi_\a} e_1=\rho_\a e_1+\kappa_\a e_3+\lambda_\a e_4,\\
A_{\xi_\a} e_2=\rho_\a e_2+\mu_\a e_3+\nu_\a e_4,\\
A_{\xi_\a} e_3=\kappa_\a e_1+\mu_\a e_2+\sigma_\a e_3,\\
A_{\xi_\a} e_4=\lambda_\a e_1+\nu_\a e_2+\sigma_\a e_4,
\end{cases}
\ee
where
\be\label{1st}
\rho_1+\sigma_1=2H\;\;\text{and}\;\;\rho_\a+\sigma_\a=0
\;\;\text{for any}\;\;2\leq\a\leq m.\nonumber
\ee
Since each shape operator has at most two distinct 
eigenvalues, each of multiplicity $2$, it follows from 
\eqref{A} that 
\be\label{munu}
\nu_\a=\pm\kappa_\a\;\;\text{and}\;\;\mu_\a=\mp\lambda_\a
\;\;\text{for any}\;\;1\leq\a\leq m.\nonumber
\ee
From this, \eqref{a1} follows directly, along with the relations: 
\be\label{a1'}
\<\a_{13},\a_{14}\>+\<\a_{23},\a_{24}\>=0,\;\,
\<\a_{13},\a_{23}\>+\<\a_{14},\a_{24}\>=0.
\ee 
Using \eqref{a1}, \eqref{a1'}, and the relation 
$\a_{11}+\a_{44}=\mathcal H/2$, we deduce that \eqref{1}, 
which now holds as an equality, is equivalent to \eqref{a3}. 
Furthermore, if \eqref{a2} is satisfied, the orthonormal four-frame 
$\{e_i\}_{1\leq i\leq4}$ satisfies all the desired properties.

Suppose now that the frame $\{e_i\}_{1\leq i\leq4}$ does 
not satisfy \eqref{a2}. Consider instead the orthonormal 
four-frame $\{\tilde e_i\}_{1\leq i\leq4}$, defined by 
$$
\tilde e_i=R_\theta e_i,\;\tilde e_j=R_\varphi e_j,\;i=1,2,\;j=3,4.
$$
Here, $R_\theta$ and $R_\varphi$ 
denote rotations on $V_1=\spa\{e_1,e_2\}$ and 
$V_2=\spa\{e_3,e_4\}$ by angles $\theta$ and 
$\varphi$, respectively. For simplicity, we set
$\tilde\a_{ij}=\a_f(\tilde e_i,\tilde e_j)$. Then, 
using \eqref{a1}, we have 
$$
\tilde\a_{ii}=\a_{ii},\;1\leq i\leq4,\;\;{\text {and}}\;\;\tilde\a_{12}=\tilde\a_{34}=0. 
$$
Since each shape operator has at most two distinct 
eigenvalues, each of multiplicity $2$, a similar 
argument as above shows that the vectors $\tilde\a_{ij}$, 
for $1\leq i,j\leq4$, satisfy \eqref{a1}, \eqref{a3} and \eqref{a1'}. 
We now claim that the angles $\theta$ and $\varphi$ 
can be chosen such that \eqref{a2} is satisfied for the frame 
$\{\tilde e_i\}_{1\leq i\leq4}$. Straightforward computations then yield 
the following relations:
\begin{align}
&\tilde\a_{14}+\tilde\a_{23}=\cos(\varphi+\theta)(\a_{14}+\a_{23})
+\sin(\varphi+\theta)(\a_{24}-\a_{13}),\label{aa1}\\
&\tilde\a_{24}-\tilde\a_{13}=-\sin(\varphi+\theta)(\a_{14}+\a_{23})
+\cos(\varphi+\theta)(\a_{24}-\a_{13}),\label{aa2}\\
&\tilde\a_{13}+\tilde\a_{24}=\cos(\varphi-\theta)(\a_{13}+\a_{24})
+\sin(\varphi-\theta)(\a_{14}-\a_{23}),\label{aa3}\\
&\tilde\a_{14}-\tilde\a_{23}=-\sin(\varphi-\theta)(\a_{13}+\a_{24})
+\cos(\varphi-\theta)(\a_{14}-\a_{23})\label{aa4}.
\end{align}
We can select angles $\sigma_1$ and $\sigma_2$ such that 
the following conditions are satisfied:
\begin{align*}
&2\cos\sigma_1\<\a_{14}+\a_{23},\a_{24}-\a_{13}\>
+\sin\sigma_2\big(\|\a_{13}-\a_{24}\|^2-\|\a_{14}+\a_{23}\|^2\big)=0,\\
&2\cos\sigma_2\<\a_{13}+\a_{24},\a_{14}-\a_{23}\>
-\sin\sigma_2\big(\|\a_{14}-\a_{23}\|^2-\|\a_{13}+\a_{24}\|^2\big)=0.
\end{align*}
By choosing 
$$
\varphi=\frac{1}{4}(\sigma_1+\sigma_2)\;\;{\text {and}}\;\;
\theta=\frac{1}{4}(\sigma_1-\sigma_2), 
$$ 
and applying \eqref{aa1}-\eqref{aa4}, 
it follows that 
$$
\<\tilde\a_{14}+\tilde\a_{23},\tilde\a_{13}-\tilde\a_{14}\>=0,
\;\,\<\tilde\a_{13}+\tilde\a_{24},\tilde\a_{14}-\tilde\a_{23}\>=0.
$$
This completes the proof of the proposition.\qed
\vspace{1.5ex}

We recall that the \emph{first normal space} of an isometric 
immersion $f\colon M^n\to\Q_c^{n+m}$ at a point 
$x\in M^n$ is defined by  
$$
N_1f(x)=\operatorname{span}\left\{\a_f(X,Y):X,Y\in T_xM\right\}.
$$

\begin{lemma}\label{nnic1}
Let $f\colon M^4\to\Q^{4+m}_c, c\geq 0$, be an 
isometric immersion of an oriented $4$-dimensional 
Riemannian manifold $M^4$. Suppose that at a point 
$x\in M^4$, there exists an oriented orthonormal 
four-frame $\{e_1,e_2,e_3,e_4\}\subset T_xM$ 
such that the second fundamental form of $f$ at $x$ 
satisfies conditions \eqref{a1}-\eqref{a3} in 
part (ii) of Proposition \ref{nnic}. 
We consider the orthonormal basis 
$\{\eta_i\}_{1\leq i\leq6}$ of the space 
of 2-vectors $\Lambda^2T_xM$, defined by the 
relations
$$
\eta_i\in\Lambda_+^2T_xM,\;\;\eta_{i+3}=\ast\eta_i, 
\;\;\text{for}\;\; 1\leq i\leq 3,
$$
and 
\begin{align*}
\eta_1=\frac{1}{\sqrt{2}}(e_{12}+e_{34}),\;
\eta_2&=\frac{1}{\sqrt{2}}(e_{13}-e_{24}),\;
\eta_3=\frac{1}{\sqrt{2}}(e_{14}+e_{23}),
\end{align*}
where $e_{ij}=e_i\wedge e_j$. 
Then the following assertions hold at $x$:
\vspace{0.5ex}

\noindent $(i)$ The matrix of the Bochner-Weitzenb\"ock 
operator $\B^{[2]}=\B_+^{[2]}\oplus\B_-^{[2]}$ at the point $x$, with respect to the 
basis $\{\eta_i\}_{1\leq i\leq6}$, is given by the direct sum 
\be\label{ssf}
\begin{bmatrix}
\mu^+_1&a^+_1&a^+_2\\
a^+_1&\mu^+_2&0\\
a^+_2&0&\mu^+_3&
\!\!\!\!\!\end{bmatrix}
\oplus
\begin{bmatrix}
\mu^-_1&a^-_1&a^-_2\\
a^-_1&\mu^-_2&0\\
a^-_2&0&\mu^-_3&
\!\!\!\!\!\end{bmatrix},\nonumber
\ee
where
\begin{align*}
\mu^\pm_1=\;&\|\a_{14}\pm\a_{23}\|^2+\|\a_{13}\mp\a_{24}\|^2,\\
\mu^\pm_2=\;&\|\a_{11}-\a_{44}\|^2+4\|\a_{13}\|^2
+2\left(\|\a_{14}\|\|\a_{23}\|\mp\<\a_{14},\a_{23}\>\right),\\
\mu^\pm_3=\;&\|\a_{11}-\a_{44}\|^2+4\|\a_{14}\|^2
+2\left(\|\a_{13}\|\|\a_{24}\|\pm\<\a_{13},\a_{24}\>\right),\\
a^\pm_1=\;&\<\a_{23}\pm\a_{14},\a_{44}-\a_{11}\>,
\;\,a^\pm_2=\<\a_{24}\mp\a_{13},\a_{44}-\a_{11}\>.
\end{align*}

\noindent $(ii)$ If $\ker \B_+^{[2]}\neq0$ at $x$, then one 
of the following conditions holds:
\begin{align}
&a_{14}+\a_{23}=0=\a_{13}-\a_{24}\;\;{\text and}\;\;a_1^+=a_2^+=0,\tag{ii1}\label{ii1}\\
&\a_{13}=\a_{24}=0,\;\a_{23}=\a_{14}\neq0\;\;{\text and}\;\;\a_{44}-\a_{11}=2\rho\,\a_{14},
\;\rho\in\R,\tag{ii2}\label{ii2}\\
&\a_{14}=\a_{23}=0,\;\a_{24}=-\a_{13}\neq0\;\;{\text and}\;\;\a_{44}-\a_{11}=2\rho\,\a_{13},
\;\rho\in\R.\tag{ii3}\label{ii3}
\end{align}

\noindent $(iii)$ If $\ker \B_-^{[2]}\neq0$ at $x$, then one 
of the following conditions holds:
\begin{align}
&\a_{14}-\a_{23}=0=\a_{13}+\a_{24}\;\;{\text and}\;\;a_1^-=a_2^-=0,\tag{iii1}\label{iii1}\\
&\a_{13}=\a_{24}=0,\;\a_{23}=-\a_{14}\neq0\;\;{\text and}\;\;a_{44}-\a_{11}=2\rho\,\a_{14},
\;\rho\in\R,\tag{iii2}\label{iii2}\\
&\a_{14}=\a_{23}=0,\;\a_{24}=\a_{23}\neq0\;\;{\text and}\;\;\a_{44}-\a_{11}=2\rho\,\a_{13},
\;\rho\in\R.\tag{iii3}\label{iii3}
\end{align}
In addition, in cases \eqref{ii2}, \eqref{ii3}, \eqref{iii2}, and \eqref{iii3}, the immersion $f$ 
has flat normal bundle at $x$, and its first normal space satisfies 
$$1\leq\dim N_1f(x)\leq2.
$$ 
\end{lemma} 
\proof 
\noindent $(i)$ By a straightforward computation 
using \eqref{WB}, the Gauss equation, and 
\eqref{a1}-\eqref{a3}, we obtain the following: 
\begin{align*}
\B^{[2]}(\eta_1)=&\;2\big(2c+2\<\a_{11},\a_{44}\>-\|\a_{13}\|^2-\|\a_{14}\|^2+
\<\a_{14},\a_{23}\>-\<\a_{13},\a_{24}\>\big)\eta_1\\
&+\<\a_{14}+\a_{23},\a_{44}-\a_{11}\>\eta_2+\<\a_{24}-\a_{13},\a_{44}-\a_{11}\>\eta_3\\
=&\;\big(\|\a_{14}+\a_{23}\|^2+\|\a_{13}-\a_{24}\|^2 \big)\eta_1 +a_1^+\eta_2+ a_2^+\eta_3,\\
\B^{[2]}(\eta_2)=&\;\<\a_{14}+\a_{23},\a_{44}-\a_{11}\>\eta_1\\
&+\big(4c+\|\a_{11}\|^2+\|\a_{44}\|^2+2\<\a_{11},\a_{44}\> -2\|\a_{14}\|^2-2\<\a_{14},\a_{23}\>\big)\eta_2\\
=&\;a_1^+\eta_1+\big(\|\a_{11}-\a_{44}\|^2+4\|\a_{13}\|^2+2\|\a_{14}\|\|\a_{23}\|
-2\<\a_{14},\a_{23}\>\big)\eta_2,\\
\B^{[2]}(\eta_3)=&\;\<\a_{24}-\a_{13},\a_{44}-\a_{11}\>\eta_1\\
&+\big(4c+\|\a_{11}\|^2+\|\a_{44}\|^2+2\<\a_{11},\a_{44}\>-2\|\a_{13}\|^2
+2\<\a_{13},\a_{24}\>\big)\eta_3\\
=&\;a_2^+\eta_1+\big(\|\a_{11}-\a_{44}\|^2+4\|\a_{14}\|^2+2\|\a_{13}\|\|\a_{24}\|
+2\<\a_{13},\a_{24}\>\big)\eta_3,
\end{align*}
and 
\begin{align*}
\B^{[2]}(\eta_4)=&\;2\big(2c+2\<\a_{11},\a_{44}\>-\|\a_{13}\|^2-\|\a_{14}\|^2-
\<\a_{14},\a_{23}\>+\<\a_{13},\a_{24}\>\big)\eta_4\\
&+\<\a_{23}-\a_{14},\a_{44}-\a_{11}\>\eta_5
+\<\a_{24}+\a_{13},\a_{44}-\a_{11}\>\eta_6\\
=&\;\big(\|\a_{14}-\a_{23}\|^2+\|\a_{13}+\a_{24}\|^2\big)\eta_4+a_1^-\eta_5+ a_2^-\eta_6,
\end{align*}
\begin{align*}
\B^{[2]}(\eta_5)=&\;\<\a_{23}-\a_{14},\a_{44}-\a_{11}\>\eta_4\\
&+\big(4c+\|\a_{11}\|^2+\|\a_{44}\|^2+2\<\a_{11},\a_{44}\>-2\|\a_{14}\|^2
+2\<\a_{14},\a_{23}\>\big)\eta_5\\
=&\; a_1^-\eta_4+\big(\|\a_{11}-\a_{44}\|^2+4\|\a_{13}\|^2+2\|\a_{14}\|\|\a_{23}\|
+2\<\a_{14},\a_{23}\>\big)\eta_5,\\
\B^{[2]}(\eta_6)=&\;\<\a_{24}+\a_{13},\a_{44}-\a_{11}\>\eta_4\\
&+\big(4c+\|\a_{11}\|^2+\|\a_{44}\|^2+2\<\a_{11},\a_{44}\>-2\|\a_{13}\|^2
-2\<\a_{13},\a_{24}\>\big)\eta_6\\
=&\; a_2^-\eta_4+\big(\|\a_{11}-\a_{44}\|^2+4\|\a_{14}\|^2+2\|\a_{13}\|\|\a_{24}\|
-2\<\a_{13},\a_{24}\>\big)\eta_6.
\end{align*}
This completes the proof of part $(i)$. 
\vspace{1ex}

$(ii)$ Using part $(i)$, we find that 
$$
\det\B_+^{[2]}=\mu^+_1\mu^+_2\mu^+_3-\mu^+_2(a_2^+)^2-\mu^+_3(a_1^+)^2.
$$
Hence, if $\ker \B_+^{[2]}\neq0$ at $x$, then
$$
\mu^+_1\mu^+_2\mu^+_3=\mu^+_2(a_2^+)^2+\mu^+_3(a_1^+)^2,
$$
or equivalently 
\be\label{zero}
\mu^+_2\big(\mu^+_3\|\a_{13}-\a_{24}\|^2-(a_2^+)^2\big)
+\mu^+_3\big(\mu^+_2\|\a_{14}+\a_{23}\|^2-(a_1^+)^2\big)
=0. 
\ee
On the other hand, by applying the Cauchy-Schwarz 
inequality, we obtain 
\begin{align*}
\mu^+_2\big(\mu^+_3\|\a_{13}-\a_{24}\|^2&-(a_2^+)^2\big)
+\mu^+_3\big(\mu^+_2\|\a_{14}+\a_{23}\|^2-(a_1^+)^2\big)\\
\geq&\;\mu^+_2\|\a_{13}-\a_{24}\|^2\big(\mu^+_3-\|\a_{44}-\a_{11}\|^2\big)\\
&+\mu^+_3\|\a_{14}+\a_{23}\|^2\big(\mu^+_2-\|\a_{44}-\a_{11}\|^2\big)\\
=&\;2\mu^+_2\|\a_{13}-\a_{24}\|^2\big(2\|\a_{14}\|^2+
\|\a_{13}\|\|\a_{24}\|+\<\a_{13},\a_{24}\>\big)\\
&+2\mu^+_3\|\a_{14}+\a_{23}\|^2\big(2\|\a_{13}\|^2+\|\a_{14}\|\|\a_{23}\|-\<\a_{14},\a_{23}\>\big).
\end{align*}
The above, together with \eqref{zero}, implies that 
\begin{align*}
&\mu^+_2\|\a_{13}-\a_{24}\|^2\big(2\|\a_{14}\|^2+
\|\a_{13}\|\|\a_{24}\|+\<\a_{13},\a_{24}\>\big)=0,\\
&\mu^+_3\|\a_{14}+\a_{23}\|^2\big(2\|\a_{13}\|^2+\|\a_{14}\|\|\a_{23}\|-\<\a_{14},\a_{23}\>\big)=0.
\end{align*}
This clearly concludes the proof 
of part $(ii)$. The proof of part $(iii)$ proceeds similarly 
and is therefore omitted.

The Ricci equation states that the normal curvature tensor  
$R^\perp$ of the submanifold $f$ is given by 
$$
R^\perp(X,Y)\xi=\a_f(X,A_\xi Y)-\a_f(A_\xi X,Y)
$$
for all $X,Y\in\mathcal X(M)$ and $\xi\in\Gamma(N_fM)$. 
We claim that 
\be\label{flatn}
\a_f(e_i,A_\xi e_j)-\a_f(A_\xi e_i,e_j)=0 \;\;\text{for all}\;\;1\leq i,j\leq4,
\ee
which clearly implies that the immersion$f$ has flat normal 
bundle at the point $x$. 

We prove \eqref{flatn} in the case where \eqref{ii2} holds, 
since the same argument applies to the other cases. 
Let $(h_{ij})$ denote the matrix of the shape  
operator $A_\xi$ with respect to the four-frame 
$\{e_1,e_2,e_3,e_4\}\subset T_xM$. 
Then it follows from part \eqref{ii2} of Lemma 
\ref{nnic1} that 
$$
h_{13}=h_{24}=0,\;h_{23}=-h_{14}\neq0\;\;\text{and}\;\;
h_{44}-h_{11}=2\rho\, h_{14}. 
$$
That \eqref{flatn} holds follows directly from these relations 
together with part \eqref{ii2} of Lemma \ref{nnic1}. 

Finally, the inequality on the dimension of the first normal 
space at $x$ follows directly from parts \eqref{ii2}, \eqref{ii3}, 
\eqref{iii2} and \eqref{iii3} of Lemma \ref{nnic1}. 
\qed
\vspace{1.5ex}

Let $f\colon M^n\to\Q_c^{n+m}$ be an isometric immersion. 
We recall that a vector $\eta$ in the normal space $N_fM(x)$ 
is called a \emph{Dupin principal normal} of the isometric 
immersion $f$ at a point $x\in M^n$ if the tangent subspace
$$
E_\eta(x)=\left\{X\in T_xM:\alpha_f(X,Y)
=\<X,Y\>\eta\;\;\text{for all}\;\,Y\in T_xM\right\}
$$
is at least two dimensional. That dimension is called the 
\emph{multiplicity} of $\eta$. The \emph{relative nullity} 
subspace $\mathcal D_f(x)$ of $f$ at a point $x\in M^n$ is the 
kernel of the second fundamental form at this point, namely
$$
\mathcal D_f(x)=\left\{X\in T_xM:\alpha_f(X,Y)
=0\;\;\text{for all}\;\,Y\in T_xM\right\}.
$$

The following result will be used in the sequel.

\begin{proposition}\label{M2compact}
Let $f\colon M^n\to\Q_c^{n+m},n\geq4,c\geq0$, 
be a substantial isometric immersion of a compact, 
simply connected, even-dimensional manifold 
with flat normal bundle. Assume that equality 
holds in \eqref{1} for $k=n/2$, and that there 
exist two distinct Dupin principal normals 
$\eta_1$ and $\eta_2$, each of multiplicity 
$k$, at every point $x\in M^n$ where 
$f$ is not umbilical, such that 
$$
T_xM=E_{\eta_1}(x)\oplus E_{\eta_2}(x).
$$
Then, $M^n$ is 
isometric to a torus $\Sf^k(r)\times \Sf^k(\sqrt{R^2-r^2})$ 
and $f$ is a composition $f=j\circ g$, where 
$g\colon M^n\to\Sf^{n+1}(R)$ is the standard 
embedding of the torus 
$\Sf^k(r)\times \Sf^k(\sqrt{R^2-r^2})$ into a 
sphere $\Sf^{n+1}(R)$, and 
$j\colon\Sf^{n+1}(R)\to\Q_c^{n+2}$ is an 
umbilical inclusion, with $R\leq1$ if $c=1$.
\end{proposition}
\proof
Consider the open subset $M^*$ of points 
where $f$ is not umbilical. 
On $M^*$, the vector fields 
$\eta_1$ and $\eta_2$ are smooth Dupin principal 
normals, with smooth corresponding distributions 
$E_1=E_{\eta_1}$ and $E_2=E_{\eta_2}$, each of rank $k$ (see 
Proposition 1 in \cite{Reck} or Lemma 10 in \cite{v3}). 

Since 
$$
2\mathcal H_f=\eta_1+\eta_2 \;\;{\text {and}}\;\, S=k(\|\eta_1\|^2+\|\eta_2\|^2),
$$ 
condition \eqref{1}, 
which now holds as equality for $k=n/2$, implies 
\be\label{eta}
\<\eta_1,\eta_2\>=-c.
\ee

The Codazzi equation for $f$ is easily seen to yield
\be\label{cod}
\<\nabla_XY,Z\>(\eta_i-\eta_j)=\<X,Y\>\nabla_Z^\perp\eta_i
\;\,\text{if}\;\,i\neq j,
\ee
for all $X,Y\in E_i,Z\in E_j,1\leq i,j\leq 2$. 
Using equation \eqref{eta} and the fact that the Dupin principal 
normal vector fields $\eta_1$ and 
$\eta_2$ are parallel in the normal connection along 
$E_1$ and $E_2$, respectively (see Proposition 1.22 in 
\cite{DT}), it follows from equation \eqref{cod} that the following 
relations hold on $M^*$:
\be\label{s}
\<\nabla_XY,Z\>(c+\|\eta_j\|^2)=0
\;\;\text{for all}\;\,X,Y\in E_i\;\;\text{and}\;\,Z\in E_j,
\;i\neq j.
\ee

Now we distinguish two cases.
\vspace{1ex}

\indent\emph{Case $c>0$}. From equation \eqref{eta}, it  
follows that $M^*=M^n$. Consequently, equation  
\eqref{s} implies that the distributions $E_1$ and $E_2$ 
are totally geodesic on $M^n$. Since $M^n$ is 
simply connected, it is well-known that $M^n$ is 
a Riemannian product $M_1^k\times M_2^k$ 
(cf.\ Theorem $8.2$ in \cite{DT}) such that 
$TM_i^k=E_i$, $i=1,2$. As the second 
fundamental form of $f$ is adapted to the 
distributions $E_1$ and $E_2$, then Theorem 
$8.4$ and Corollary $8.6$ in \cite{DT} imply 
that the submanifold is an extrinsic product 
of isometric immersions, each of which is totally 
umbilical. Therefore, the submanifold is a 
torus $\Sf^k(r)\times\Sf^k(\sqrt{R^2-r^2})$ in a 
sphere $\Sf^{n+1}(R)\subset\Sf^{n+2}$.
\vspace{1ex}

\indent\emph{Case $c=0$}. Consider the open subset
$$
M_+=\left\{x\in M^*: \eta_1(x)\neq0 \,\;\text{and}\,\;\eta_2(x)\neq0\right\}.
$$ 

Suppose that $M_+$ is nonempty. From \eqref{s}, 
it follows that that the distributions 
$E_1$ and $E_2$ are totally geodesic on $M_+$. 
It is well-known that $M_+$ is locally a Riemannian 
product $M_1^k\times M_2^k$ (cf.\ Theorem $8.2$ 
in \cite{DT}) such that $TM_i^k=E_i$, $i=1,2$. 
Since the second fundamental form of $f$ is adapted 
to the distributions $E_1$ and $E_2$, Theorem $8.4$ 
and Corollary $8.6$ in \cite{DT} imply that the submanifold 
is locally an extrinsic product of isometric immersions, 
each of which is totally umbilical. Consequently, 
the submanifold is locally, on $M_+$, a torus 
$\Sf^k(r)\times\Sf^k(\sqrt{R^2-r^2})$ in a sphere 
$\Sf^{n+1}(R)\subset \R^{n+2}$. 
\vspace{1ex}

Now, suppose that the interior $\operatorname{int}(M_0)$
of the subset $M_0=M^*\smallsetminus M_+$ 
is nonempty. Assume that $\eta_2=0$ 
on a connected component $U$ of $\operatorname{int}(M_0)$. 
Then $f$ has constant index of relative nullity $k$ 
on $U$. From \eqref{s}, it follows that the relative 
nullity distribution $\mathcal D_f=E_2$ is totally geodesic. 
Let $C_T\colon\mathcal D_f^\perp\to\mathcal D_f^\perp$ 
be the associated splitting tensor for any 
$T\in\mathcal D_f$ (see \cite[p. 186]{DT}). Since 
$$
\a_f(X,Y)=
\<X,Y\>\eta_1\;\;{\text{for all}}\;\,X\in\mathcal X(M^n),Y\in\mathcal D_f^\perp,
$$
the Codazzi equation
$$
(\n_X^\perp\a_f)(Y,T)=
(\n_T^\perp\a_f)(X,Y)\;\;{\text{for all}}\;\,X,Y\in\mathcal D_f^\perp, T\in\mathcal D_f
$$
implies that 
$$
C_T=\<\n\log\|\eta_1\|, T\>Id,
$$ 
where $Id$ is the identity map on the conullity distribution 
$\mathcal D_f^\perp$. This shows that the conullity 
distribution is umbilical, and consequently, 
integrable. 

Let $\varSigma^k$ be a leaf of $\mathcal D_f^\perp$. By 
Proposition 7.6 in \cite{DT}, $f$ is locally a 
generalized cone over an isometric immersion 
$g\colon\varSigma^k\to \Q_{\tilde c}^{k+m}$ such 
that $f\circ j_{\varSigma}=i\circ g$, where 
$j_{\varSigma}\colon\varSigma^k\to M^n$ is the inclusion and 
$i\colon\Q_{\tilde c}^{k+m}\to\R^{n+m}$ is an 
umbilical inclusion. The second fundamental 
form of the immersion $f\circ j$ is given by 
$$
\a_{f\circ j}(X,Y)=
\<X,T\>\big(\eta_1+f_*\n\log\|\eta_1\|\big)\;\;{\text{for all}}\;\,X,Y\in\mathcal X(\varSigma^k).
$$
Thus, $f\circ j_{\varSigma}$ is umbilical and $g(\varSigma^k)$ is 
a sphere $\Sf^k(r)$ centered at a point 
$x_0\in\R^{n+m}$ with radius 
$$
r=\frac{1}{\sqrt{\|\eta_1\|^2+\|\n\log\|\eta_1\|\|^2}}. 
$$ 
If the umbilical submanifold $\Q_{\tilde c}^{k+m}$ 
is totally geodesic in $\R^{n+m}$, then the 
submanifold $f$ is a $k$-cylinder over the 
sphere $\Sf^k(r)$. Alternatively, if 
$\Q_{\tilde c}^{k+m}$ is a sphere centered 
at a point $\tilde x_0\in\R^{n+m}$, then 
$f$ is a $(k-1)$-cylinder over a submanifold 
$N^{k+1}$, which is a cone over the sphere 
$\Sf^k(r)$ with its vertex at $\tilde x_0\neq x_0$.

Consequently, on $\operatorname{int}(M_0)$ the 
submanifold is locally a $k$-cylinder 
in $\R^{n+1}$ over a sphere $\Sf^k(r)$, or 
a $(k-1)$-cylinder over a submanifold 
$N^{k+1}$ in $\R^{k+2}$, which is a cone 
over a sphere $\Sf^k(r)$. 

In the later case, $f$ is given locally on 
$\operatorname{int}(M_0)$ by
$$
f(x,w)=j(x)+w, \;\; (x,w)\in N_i\Sf^{k+m}(R),
$$
where the normal bundle $N_i\Sf^{k+m}(R)$ 
of the inclusion $i\colon\Sf^{k+m}(R)\to\R^{n+m}$ 
is regarded as subbundle of $N_{i\circ j}\Sf^k(r)$. 
Here, $j\colon\Sf^k(r)\to\Sf^{k+m}(R)$ is an 
umbilical inclusion with $r<R$. Equivalently, 
$f$ is locally parametrized by
$$
f(x,t_0,t_1,\dots,t_{k-1})=t_ 0j(x)+
\sum_{i=1}^{k-1}t_\ell v_\ell,\;\,t_0>0,\; x\in\Sf^k(r),
$$
where 
$$
\Sf^{k+m}(R)\subset\R^{k+m+1}, \;\;
\R^{n+m}=\R^{k+m+1}\oplus\R^{k-1},
$$ 
and 
$\{v_\ell\}_{1\leq\ell\leq k-1}$ is an orthonormal 
basis of $\R^{k-1}$. The Laplacian operator 
$\Delta_M$ of $M^n$ is given by
$$
\Delta_M=\frac{k}{t_0}\frac{\partial}{\partial t_0}+
\frac{\partial^2}{\partial t_0^2}+\frac{1}{t^2_0}\,\Delta_{\Sf^k(r)}
+\sum_{\ell=1}^{k-1}\frac{\partial^2}{\partial t^2_\ell}.
$$
Since $S=a/t^2_0$, where $a$ is a positive 
constant, it follows that 
$$
\Delta_M S=-\frac{a(n-6)}{t_0^4}. 
$$
Consequently, either $\Delta_M S\geq0$ 
or $\Delta_M S\leq0$ on $\operatorname{int}(M_0)$, 
depending on the dimension $n$. 

Since $S$ is locally constant on $M_+$, 
continuity implies that $S$ is either 
superharmonic or subharmonic 
function on $M^n$. By the maximum 
principle, it follows that $S$ is a positive 
constant on $M^n$. Consequently, $M^*=M^n$, 
and therefore $M^n=M_+\cup M_0$. 

Thus, on $M_0$, one of the Dupin principal 
normal vector fields vanishes, while the 
other has constant length $S/k$. On the 
other hand, both principal normal fields have 
constant and positive length on each 
connected component of $M_+$.
By continuity and the compactness of 
$M^n$, it follows that $M_+=M^n$. 
Consequently, using a similar argument 
as in \emph{Case} $c>0$, we conclude that $M^n$ 
is isometric to a torus 
$\Sf^k(r)\times \Sf^k(\sqrt{R^2-r^2})$, and  
the immersion $f$ is the standard embedding 
into a sphere $\Sf^{n+1}(R)\subset \R^{n+2}$.
\qed
\vspace{1.5ex}

\noindent\emph{Proof of Theorem \ref{k=2}:}
Part $(i)$ of Proposition \ref{nnic} implies that $M^4$ 
has nonnegative isotropic curvature and that the 
Bochner-Weitzenb\"ock operator $\B^{[2]}$ is 
nonnegative at every point. The proof of the 
theorem is divided into three cases.
\vspace{1ex}

\indent\emph{Case I}. 
We begin by proving the theorem for simply 
connected submanifolds. Since $M^4$ has nonnegative 
isotropic curvature, Theorem 4.10 in \cite{MW} implies 
that one of the following holds: 
\begin{enumerate}[topsep=1pt,itemsep=1pt,partopsep=1ex,parsep=0.5ex,leftmargin=*, label=(\roman*), align=left, labelsep=-0.5em]
\item[(a)] $M^4$ carries a metric of positive isotropic 
curvature. 
\item[(b)] $M^4$ is diffeomorphic to a product 
$\Sf^2\times\varSigma^2$, where $\varSigma^2$ is a 
compact surface.
\item[(c)] $M^4$ is a K\"ahler manifold biholomorphic 
to $\CP^2$.
\end{enumerate} 

We analyse each case separately as follows.
\vspace{1ex}

\noindent\emph{Case} (a). We assert 
that the manifold $M^4$ is diffeomorphic to $\Sf^4$. 
Given that $M^4$ is simply connected by assumption, 
the main result in \cite{MM} implies that $M^4$ is 
homeomorphic to $\Sf^4$. Furthermore, $M^4$ 
is locally irreducible. If this were not the case, then 
Theorem 3.1 in \cite{MW} would imply that $M^4$ is 
isometric to a Riemannian product of two compact 
surfaces, which leads to a contradiction. Hence, 
$M^4$ is locally irreducible, and by Theorem 
2 in \cite{BSacta}, one of the following cases must hold:
\begin{enumerate}[topsep=1pt,itemsep=1pt,partopsep=1ex,parsep=0.5ex,leftmargin=*, label=(\roman*), align=left, labelsep=-0.5em]
\item[(i)] $M^4$ is diffeomorphic to a spherical space form. 
\item[(ii)] $M^4$ is a K\"ahler manifold biholomorphic to 
$\CP^2$.
\item[(iii)] $\;M^4$ is isometric to a compact symmetric space.
\end{enumerate} 

Since $M^4$ is homeomorphic to $\Sf^4$, case (ii) above 
is excluded. In case (i), the manifold is clearly diffeomorphic 
to $\Sf^4$. Moreover, the only  
compact symmetric spaces of dimension four are the round sphere, the 
product of two 2-dimensional spheres, the complex projective plane, and 
the quaternionic projective line $\mathbb HP^1\cong \Sf^4$. 
Consequently, in case (iii), the manifold $M^4$ must also 
be diffeomorphic to $\Sf^4$. 
\vspace{1ex}

\noindent\emph{Case} (b). Since $M^4$ is 
simply connected, the surface $\varSigma^2$ must be 
diffeomorphic to $\Sf^2$, which implies that $M^4$ is 
diffeomorphic to $\Sf^2\times\Sf^2$. Consequently, 
$\beta_\pm(M)=1$. It then follows from part $(ii)$ of 
Proposition \ref{nnic} that equality holds in \eqref{1}. 
Moreover, at each point $x\in M^4$, there exists an 
oriented orthonormal four-frame 
$\{e_1,e_2,e_3,e_4\}\subset T_xM$ such that 
the second fundamental form $\a_f$ of $f$ 
satisfies conditions \eqref{a1}-\eqref{a3}. 

In addition, there exist a nontrivial 
self-dual harmonic 2-form $\omega_+$ and 
a nontrivial anti-self-dual harmonic 2-form $\omega_-$. 
By the Bochner-Weitzenb\"ock 
formula, both $\omega_+$ and $\omega_-$
are parallel and 
$$
\<\B^{[2]}(\omega_\pm),\omega_\pm\>=0
$$ 
everywhere. Consequently, $\ker \B_+^{[2]}\neq0$ 
and $\ker \B_-^{[2]}\neq0$ at every point. Hence, both 
parts $(ii)$ and $(iii)$ of Lemma \ref{nnic1} hold everywhere.

We claim that the immersion $f$ has flat normal bundle at 
every point $x \in M^4$. Suppose first that both cases 
\eqref{ii1} and \eqref{iii1} hold simultaneously at $x$. 
Then, it follows from \eqref{a1} and 
\eqref{a3} that $c=0$, and hence $f$ is totally geodesic 
at this point. If, instead, either \eqref{ii1} or \eqref{iii1} 
does not hold at $x$, then Lemma~\ref{nnic1} implies 
that $f$ has flat normal bundle at this point. 
This proves the claim that $f$ has flat normal bundle.

Since by part $(ii)$ of Proposition \ref{nnic},  
the shape operator, associated with any normal direction, 
has at most two distinct eigenvalues, each of multiplicity 
$2$ at every point. This together with the flatness of the 
normal bundle implies that there exist two distinct Dupin 
principal normals, $\eta_1$ and $\eta_2$ each of multiplicity 
$2$, at every point $x\in M^4$ where $f$ is not umbilical, 
such that 
$$
T_xM=E_{\eta_1}(x)\oplus E_{\eta_2}(x).
$$

It then follows from Proposition \ref{M2compact} that the 
submanifold is as described in part $(iii$-$a)$ of the 
theorem. 
\vspace{1ex}

\noindent\emph{Case} (c). In this case, either 
$$
\beta_+(M)=1\;\;\text{and}\;\;\beta_-(M)=0, \;\;\text{or}\;\; 
\beta_+(M)=0\;\;\text{and}\;\;\beta_-(M)=1.
$$ 
We will treat only the former case, as 
the latter can be handled in a similar 
manner. It then follows from part $(ii)$ of 
Proposition \ref{nnic} that equality holds in \eqref{1}. 
Moreover, at each point $x\in M^4$, there exists an 
oriented orthonormal four-frame 
$\{e^x_1,e^x_2,e^x_3,e^x_4\}\subset T_xM$ such that 
the second fundamental form $\a_f$ of $f$ 
satisfies conditions \eqref{a1}-\eqref{a3}. 

Clearly, there exists a nontrivial self-dual harmonic 2-form 
$\omega_+$. By the Bochner-Weitzenb\"ock 
formula, the form $\omega_+$ is parallel, 
and 
$$
\<\B^{[2]}(\omega_+),\omega_+\>=0
$$ 
at every point. Hence, $\ker \B_+^{[2]}\neq0$ at every point, 
and part $(ii)$ of Lemma \ref{nnic1} applies.

We consider the subset 
$$
M_1=\left\{x\in M^4: \mu^+_1(x)=0\;\;\text{and}\;\;\mu^+_2(x)\mu^+_3(x)>0\right\},
$$
where $\mu^+_1,\mu^+_2$, and $\mu^+_3$ are as given in Lemma \ref{nnic1}. 
It follows that
$$
M_1=\left\{x\in M^4: \mu^+_1(x)=0\;\;\text{and}\;\;\|\a_{11}-\a_{44}\|^2+4\|\a_{13}\|^2
+4\|\a_{14}\|^2>0\right\}.
$$

We claim that this subset is nonempty.
Suppose, to the contrary, that 
$M_1$ is empty. Then, at every point of $M^4$, one 
of the following conditions must hold:
\be\label{n0}
\a_{14}+\a_{23}\neq0, \;\text {or}\;\a_{13}-\a_{24}\neq0,  
\;\text {or}\;\|\a_{11}-\a_{44}\|^2+4\|\a_{13}\|^2
+4\|\a_{14}\|^2=0.
\ee
We claim that the immersion $f$ has flat normal bundle. Indeed, 
if the third condition in \eqref{n0} holds at a point $x$, then 
\eqref{a1} and \eqref{a2} imply that $c=0$, and hence $f$ 
is totally geodesic at $x$. If, instead, either the first 
or the second condition in \eqref{n0} holds at a point $x$, 
then part $(ii)$ of Lemma \ref{nnic1} shows that either case 
\eqref{ii2} or case \eqref{iii2} holds at that point. In either 
situation, Lemma \ref{nnic1} further implies that the 
submanifold $f$ has flat normal bundle and satisfies 
$$
\dim N_1f\leq2
$$ 
at all points of $M^4$. 

By part $(ii)$ of Proposition \ref{nnic},  
the shape operator, associated with any normal direction, 
has at most two distinct eigenvalues, each of multiplicity 
$2$ at every point. Combined with the flatness of the 
normal bundle, this implies the existence of two distinct Dupin 
principal normals, $\eta_1$ and $\eta_2$, each of multiplicity 
$2$, at every point $x\in M$ where $f$ is not umbilical, 
such that 
$$
T_xM=E_{\eta_1}(x)\oplus E_{\eta_2}(x).
$$
Since the submanifold is substantial, Proposition 
\ref{M2compact}  implies that its codimension is $m=2$. 
Consequently, the complex projective plane $\CP^2$ would admit 
an immersion into $\R^6$ or in $\Sf^6$. This, however, 
contradicts known topological obstructions 
(see, for instance, \cite[pp. 290-291]{eft}). 
Hence, $M_1$ must be nonempty.

From part $(i)$ of Lemma \ref{nnic1}, 
it follows that at every point of $M_1$ 
we have  
$$
\mu_1^+=0, \;\; 
\mu_2^+=\mu_3^+>0.
$$ 
Consequently, the kernel of $\B_+^{[2]}$ at each 
point $x\in M_1$ is spanned by the vector 
$e^x_1\wedge e^x_2+e^x_3\wedge e^x_4$.

Let $Z$ be the dual to the self-dual 
form $\omega_+$. Since this form is parallel, 
we may normalize it (after multiplying by a 
constant if necessary) so that $\|\omega_+\|=\sqrt{2}$.
Fix a point $x_0\in M_1$. At this point, we may 
assume that  
$$
Z_{x_0}=e^{x_0}_1\wedge e^{x_0}_2+e^{x_0}_3\wedge e^{x_0}_4. 
$$
Now, consider the almost complex structure 
$J_{x_0}\colon T_{x_0}M\to T_{x_0}M$ given by
$J_{x_0}e^{x_0}_1=e^{x_0}_2$ and $J_{x_0}e^{x_0}_3=e^{x_0}_4$. 
Then, we have
$$
\omega_+(v,w)=\<J_{x_0}v,w\>\;\;\text {for all}\;\, v,w\in T_{x_0}M. 
$$

Moreover, we define the skew-symmetric endomorphism $J$ 
of the tangent bundle of $M^4$ such that 
$$
\omega_+(X,Y)=\<JX,Y\>\;\;\text {for all}\;\, X,Y\in\mathcal X(M). 
$$
Clearly, $J$ is parallel because $\omega_+$ is parallel. 
We claim that $J$ is orthogonal; that is, $\|J_xv\|=\|v\|$ 
for every point $x\in M^4$ and every $v\in T_xM$. Indeed, let 
$V$ be a parallel vector field along a curve $c\colon[0,1]\to M$ 
such that $c(0)=x,c(1)=x_0$ and $V(0)=v$. Obviously, 
the vector field $W=JV$ is also parallel along $c$. Using 
that $J_{x_0}$ is orthogonal we have 
$$
\|J_xv\|=\|W(0)\|=\|W(1)\|=\|J_{x_0}V(1)\|=\|V(1)\|=\|V(0)\|=\|v\|,
$$
which proves the claim. Since $J$ is both 
skew-symmetric and orthogonal, it defines an 
almost complex structure that is also parallel. Hence, the 
triple $(M^4,\<\cdot,\cdot\>,J)$ is a K\"ahler manifold. 

We now claim that, at every point $x\in M^4$, we have
\be\label{Jx}
\a_f(J_xv,J_xw)=\a_f(v,w)\;\;\text {for all}\;\; v,w\in T_xM.
\ee

To prove this, we distinguish two cases.

First, let $x$ be an arbitrary point in $M_1$. Since  
the kernel of $\B_+^{[2]}$ at $x$ is spanned by the vector 
$e^x_1\wedge e^x_2+e^x_3\wedge e^x_4$, and since  
$\|Z_x\|=\sqrt{2}$, it follows that 
$$
Z_x=\pm\left(e^x_1\wedge e^x_2+e^x_3\wedge e^x_4\right).
$$ 
Because $Z$ is the dual of the self-dual 
form $\omega_+$, we have 
$$
\<J_xv,w\>=\<\<Z_x,v\wedge w\>\>\;\;\text {for all}\;\, v,w\in T_xM. 
$$
Consequently,  
$$
J_xe^x_1=\pm e^x_2 \;\;\text {and}\;\; J_xe^x_3=\pm e^x_4
$$ 
at each point $x\in M_1$. On the other hand, we have  
\be\label{f}
\a_{14}+\a_{23}=\a_{13}-\a_{24}=0
\ee
at every point of $M_1$. Combining \eqref{a1} with \eqref{f}, 
we obtain 
\be\label{Jij}
\a_f(J_xe^x_i,J_xe^x_j)=\a_{ij}\;\;\text {for all}\;\; 1\leq i,j\leq4, 
\ee
and thus \eqref{Jx} follows from \eqref{Jij} by linearity.

Having analyzed the case $x\in M_1$, we now 
turn to the complementary case $x\notin M_1$.
In this case, \eqref{n0} holds at $x$.

If the third condition in \eqref{n0} is satisfied, then \eqref{a1} 
and \eqref{a3} imply that $c=0$, and hence $f$ is totally 
geodesic at $x$. In this situation, \eqref{Jx} holds trivially.

Next, suppose that the first or the second condition in 
\eqref{n0} holds at $x$. We treat only the first case, since 
the second can be handled in a similar way, and 
will therefore be omitted. Then part $(ii)$ of Lemma \ref{nnic1} 
shows that case \eqref{ii2} holds at that point. Together with 
\eqref{ii2}, it follows from 
part $(i)$ of this lemma that the matrix of 
the Bochner-Weitzenb\"ock operator $\B_+^{[2]}$ at 
$x$, with respect to the basis $\{\eta_i\}_{1\leq i\leq3}$ 
of $\Lambda_+^2T_xM$, is given by  
\be\label{ssf}
4\|\a_{14}\|^2\begin{bmatrix}
1&\rho&0\\
\rho&\rho^2&0\\
0&0&\rho^2+1&
\!\!\!\!\!\end{bmatrix},\nonumber
\ee
where $\rho$ is the constant appearing in case \eqref{ii2} of  
Lemma \ref{nnic1}. Therefore, the kernel of $\B_+^{[2]}$ 
at $x$ is spanned by the vector 
$-\rho\eta_1+\eta_2$, and since 
$\|Z_x\|=\sqrt{2}$, it follows that 
$$
Z_x=\lambda\left(-\rho\eta_1+\eta_2\right),
$$ 
where the number $\lambda$ satisfies
\be\label{lam}
\lambda^2=\frac{2}{\rho^2+1}. 
\ee
Because $Z$ is the dual of the self-dual 
form $\omega_+$, we have 
$$
\<J_xv,w\>=\<\<Z_x,v\wedge w\>\>\;\;\text {for all}\;\, v,w\in T_xM. 
$$
Using this relation together with \eqref{lam}, we obtain  
\begin{align*}
J_xe^x_1&=\frac{\lambda}{\sqrt{2}} \left(-\rho\,e^x_2+e^x_3\right), 
\;\; J_xe^x_2=\frac{\lambda}{\sqrt{2}} \left(\rho\,e^x_1-e^x_4\right),\\
J_xe^x_3&=-\frac{\lambda}{\sqrt{2}} \left(e^x_1+\rho\,e^x_4\right), 
\;\; J_xe^x_4=\frac{\lambda}{\sqrt{2}} \left(e^x_2+\rho\,e^x_3\right).
\end{align*}
From these expressions, and combining \eqref{a1} with \eqref{ii2}, 
we obtain  \eqref{Jij}. Hence \eqref{Jx} follows directly.

Therefore, the second fundamental form of the submanifold 
satisfies 
\be\label{f1}
\a_f(JX,JY)=\a_f(X,Y)\;\;\text{for all}\;\,X,Y\in\mathcal X(M).
\ee
Immersions satisfying condition \eqref{f1} have 
parallel second fundamental form 
(see Theorem 4 in \cite{f2}). Moreover, under the 
immersion $f$, each geodesic of $M^4$ is 
mapped into a plane circle. Such submanifolds 
were classified in \cite{f2} and \cite{Sk}. From 
this classification, it follows that the submanifold 
is as described in part $(iii$-$b)$ of the theorem, 
thereby completing the proof for the case 
where the manifold is simply connected.
\vspace{1ex} 

\indent\emph{Case II}. 
Now suppose that the fundamental group of $M^4$ is finite.
Then we claim that $M^4$ is simply connected. Consider the 
universal covering $\pi\colon\tilde M^4\to M^4$. Since the 
fundamental group of $M^4$ is finite, $\tilde M^4$ must 
be compact. Moreover, the isometric immersion 
$\tilde f=f\circ\pi$ satisfies \eqref{1}. Therefore, from 
\emph{Case I}, we conclude that either $\tilde M^4$ is 
diffeomorphic to $\Sf^4$, or the submanifold $\tilde f$ 
is as described in parts $(i)$ or $(ii)$ of the theorem. 

Assume first that $\tilde M^4$ is diffeomorphic to $\Sf^4$. 
Arguing as in the proof of \emph{Case I}, we conclude 
that $M^4$ is locally irreducible. Then, by Theorem 2
in \cite{BSacta}, $M^4$ is either diffeomorphic to a 
spherical space form or isometric to a compact 
symmetric space. In the former case, it 
is clear that the manifold $M^4$ must be diffeomorphic 
to $\Sf^4$. In the latter case, since the only 4-dimensional 
compact symmetric spaces are the round spheres, the 
product of two 2-dimensional spheres, the complex plane, and  
the quaternionic projective line $\mathbb HP^1$, 
it follows that $M^4$ must also be diffeomorphic to $\Sf^4$. 

Now, suppose that the submanifold $\tilde f$ is as described 
in parts $(i)$ or $(ii)$ of the theorem.
In either case, the covering map $\pi\colon\tilde M^4\to M^4$ 
must be a diffeomorphism. This observation completes the 
proof of the theorem in this case.
\vspace{1ex}

\indent\emph{Case III}. Finally, suppose that the fundamental 
group $\pi_1(M)$ of $M^4$ is infinite. We first claim 
that $M^4$ is locally reducible. If this were not the 
case, then by Theorem 2 in \cite{BSacta}, it would 
follow that either $M^4$ is diffeomorphic to a spherical space 
form, or its universal cover $\tilde M^4$ is biholomorphic to 
$\CP^2$, or isometric to a compact symmetric space. 
Each of these possibilities contradicts the fact that 
$\tilde M^4$ is not compact. Therefore, $M^4$ must be 
locally reducible. It then follows from Theorem 3.1 in 
\cite{MW} that the universal cover $\tilde M^4$ is 
isometric to one of the following Riemannian products: 
\begin{enumerate}[topsep=1pt,itemsep=1pt,partopsep=
1ex,parsep=0.5ex,leftmargin=*, label=(\roman*), align=left, labelsep=-0.5em]
\item[(i)] $(\R^{n_0},g_0)\times(N_1^{n_1},g_1)\times(N_2^{n_2},g_2)$, 
where $n_0\geq0, g_0$ is the flat Euclidean metric, and either 
$n_i=2$ and $N_i=\Sf^2, i=1,2$, has nonnegative Gaussian curvature, 
or else $n_i=3$ and $N_i$ is compact with nonnegative Ricci 
curvature. 
\item[(ii)] $(\varSigma^2,g_\varSigma)\times(N^2,g_N)$, 
where $\varSigma^2$ is a surface whose Gaussian curvature 
$K_\varSigma$ is negative at some point, and $N^2$ is a compact 
surface with positive Gaussian curvature $K_N$. 
\end{enumerate}

First, we claim that case (ii) above cannot occur. Suppose, 
for the sake of contradiction, that it does. Clearly the 
immersion $\tilde f=f\circ\pi$ satisfies \eqref{1}, where 
$\pi\colon\tilde M^4\to M^4$ is the covering map. 
We consider the isometric immersions 
$$
f_\varSigma=\tilde f\circ \mathfrak i_\varSigma\;\; \text{and}\;\;
f_N=\tilde f\circ \mathfrak i_N,
$$ 
where $\mathfrak{i}_\varSigma\colon\varSigma\to\tilde M$ and 
$\mathfrak {i}_N\colon N\to\tilde M$ are the inclusion maps. Since both 
immersions $\mathfrak i_\varSigma, \mathfrak i_N$ are totally geodesic, it follows 
that the second fundamental form of $f_\varSigma$ and $f_N$ 
are given, respectively, by 
\begin{align}
\a_{f_\varSigma}(X,Y)&=\a_{\tilde f}(\mathfrak i_{\varSigma_*} X,\mathfrak i_{\varSigma_*}Y),\;X,Y\in T\varSigma,\nonumber\\
\a_{f_N}(V,W)&=\a_{\tilde f}(\mathfrak i_{N_*}V,\mathfrak i_{N_*}W),\;V,W\in TN. \nonumber
\end{align}
Given that $\tilde M^4$ is isometric to the Riemannian product 
$(\varSigma^2,g_\varSigma)\times(N^2,g_N)$, it follows from 
the above that 
\begin{align}
S_{\tilde f}&=S_{f_\varSigma}+S_{f_N}+
2\sum_{i=1,2}\sum_{j=3,4}\|\a_{\tilde f}(\mathfrak i_{\varSigma_*} e_i,\mathfrak i_{N_*}e_j)\|^2\label{Stilde}
\end{align}
and 
\begin{align}
4\tilde H^2&=H^2_{f_\varSigma}+H^2_{f_N}
+2\<\mathcal H_{f_\varSigma},\mathcal H_{f_N}\>, \label{Htilde}
\end{align}
where $\{e_i\}_{1\leq i\leq4}$ is an orthonormal frame such that 
$e_1,e_2\in T\varSigma$, $e_3,e_4\in TN$, and 
$\tilde H$ is the mean 
curvature of the immersion $\tilde f$. 
The squared lengths of 
the traceless parts $\Phi_{f_\varSigma}$ and $\Phi_{f_N}$ of the 
second fundamental forms of the immersions $f_\varSigma$ 
and $f_N$ are given by 
$$
\|\Phi_{f_\varSigma}\|^2=S_{f_\varSigma}-2H^2_{f_\varSigma}\;\;
{\text {and}}\;\; \|\Phi_{f_N}\|^2=S_{f_N}-2H^2_{f_N}, 
$$
respectively. Taking into account \eqref{Stilde}, \eqref{Htilde} 
and the above, inequality \eqref{1} for $\tilde f$ can be 
equivalently written as 
\be\label{ph1}
\|\Phi_{f_\varSigma}\|^2+ \|\Phi_{f_N}\|^2
+2\sum_{i=1,2}\sum_{j=3,4}\|\a_{\tilde f}(\mathfrak i_{\varSigma_* }e_i,\mathfrak i_{N_*}e_j)\|^2
\leq4\big(c+\<\mathcal H_{f_\varSigma},\mathcal H_{f_N}\>\big). 
\ee

On the other hand, from the Gauss equation and \eqref{1}, 
it follows that the scalar curvature $\tau$ of $\tilde M^4$ satisfies 
$$
\tau\geq8(c+\tilde H^2).
$$ 
Since $\tau=2(K_\varSigma+K_N)$, 
we then obtain the inequality 
\be\label{sc}
K_\varSigma+K_N\geq4(c+\tilde H^2).
\ee

Using \eqref{Htilde} and the Gauss equation for both $f_\varSigma$ 
and $f_N$, we find that \eqref{sc} is written as
\be\label{ph2}
\|\Phi_{f_\varSigma}\|^2+ \|\Phi_{f_N}\|^2
+4\big(c+\<\mathcal H_{f_\varSigma},\mathcal H_{f_N}\>\big)\leq0.
\ee
Then, from \eqref{ph1} and \eqref{ph2} we obtain 
$\Phi_{f_\varSigma}=0$. Hence, the immersion $f_\varSigma$ 
it totally umbilical, and consequently $K_\varSigma\geq0$, 
which is clearly a contradiction. 

Thus, only case (i) holds. Clearly, $n_0\neq3$. We now claim 
that $n_0=1$. If $n_0=4$, then $M^4$ is flat. By the Gauss 
equation, we obtain  
$$
S=12c+16H^2.
$$ 
Together with \eqref{1}, this implies that $c=0$ and 
$f$ is minimal, which leads to a contradiction. 
If $n_0=0$ or $n_0=2$, then the curvature operator 
of $M^4$ is nonnegative. 
In this case, by Theorem 1.3 in \cite{H}, $M^4$ is 
diffeomorphic to a quotient of one of the spaces $\Sf^4$, 
$\CP^4$, $\R^1\times\Sf^3$, $\Sf^2\times\Sf^2$, 
$\R^2\times\Sf^2$, or $\R^4$ by a finite group of 
fixed-point-free isometries in the standard metric. 
This contradicts the fact that $\pi_1(M)$ is infinite. 
Therefore, $n_0=1$ and $\tilde M$ splits isometrically 
as $\R\times N^3$, where $N^3$ is compact, simply 
connected, and has nonnegative Ricci curvature. 
From Theorem 1.2 in \cite{H}, we then conclude 
that $N^3$ is diffeomorphic to $\Sf^3$.\qed
\vspace{1.5ex}

\noindent\emph{Proof of Theorem \ref{ck=2}:} 
Let $n\geq5$, and suppose that the submanifold $f$ is 
not as described in parts $(i)$-$(iii)$ of Theorem 
\ref{th1}. Then, the homology groups of $M^n$ satisfy 
$$
H_p(M^n;\Z)=0\;\,\text{for all}\;\,2\leq p\leq n-2
\;\;\text{and}\;\, H_1(M^n;\Z)\cong\Z^{\beta_1(M)}. 
$$
Since $M^n$ is oriented, we have $H^n(M^n;\Z)\cong\Z$. 
Applying the universal coefficient theorem, we deduce that
$H_{n-1}(M^n;\Z)$ is torsion-free. By Poincar\'e 
duality, it follows that 
$H_{n-1}(M^n;\Z)\cong\Z^{\beta_1(M)}$. 

In the case $n=4$, the result follows directly 
from Theorem \ref{k=2}.\qed 

\subsection{Examples of submanifolds satisfying condition \eqref{1}}

We now describe a method for constructing geometrically 
distinct submanifolds that are diffeomorphic either to the 
sphere $\Sf^n$ or to the torus $\Sf^1 \times \Sf^{n-1}$, 
and which, in addition, satisfy condition \eqref{1} in 
the case $c=0$.

\begin{proposition}\label{ex}
Let $g\colon N^{n-1}\to\R^{m_1}, n\geq4$, be 
an isometric immersion of a manifold $N^{n-1}$ satisfying 
condition \eqref{1} for an integer $k=\ell-1$, with 
$2\leq\ell\leq n-2$, as a strict inequality at every point; that is,  
$$
S_g<a(n-1,\ell-1,H_g,0). 
$$
Consider any closed unit-speed curve 
$\gamma\colon\Sf^1\to\R^{m_2}$ whose first curvature 
$\kappa_1$ satisfies 
\be\label{k1}
\kappa_1^2\leq
\frac{n-\ell}{n-\ell-1}\Big(\min\big(a(n-1,\ell-1,H_g,0)-S_g\big)\Big).
\ee
Then, the product immersion 
$$
f=\gamma\times g\colon\Sf^1\times N^{n-1}\to\R^{m_1+m_2}
$$ 
satisfies condition \eqref{1} for $k=\ell$. 
\end{proposition}
\proof
The squared length $S_f$ of the second fundamental form, 
and the mean curvature $H_f$ of the product immersion 
$f$ are given by 
$$
S_f=\kappa_1^2+S_g\;\;{\text{and}}\;\,
n^2H^2_f=\kappa_1^2+(n-1)^2H^2_g. 
$$
From these expressions, it follows that condition 
\eqref{1} for $f$, with $k=\ell$ and $c=0$, is equivalent 
to inequality \eqref{k1}, thereby completing the proof.\qed
\vspace{1.5ex}

Next, we show that there exist many isometric 
immersions $g \colon N^{n-1} \to \R^m$ satisfying 
the hypotheses of Proposition \ref{ex}.
In particular, we construct geometrically distinct 
immersions of $\Sf^n$ into the Euclidean space 
$\R^{n+1}$ that satisfy condition \eqref{1} for all $n\geq 3$ 
and $1 \leq k \leq n-1$, either with strict inequality 
or with equality.

\begin{proposition}\label{ell}
Let $f\colon M^n\to\R^{n+1},n\geq3$, be an ovaloid 
in the Euclidean space with principal curvatures 
$0<\lambda_1\leq\dots\leq\lambda_n$. If 
\be\label{lambda}
\max\lambda _n \leq\min\lambda _1\Big({\frac{n}{n-k}}\Big)^{1/2}
\ee
for some integer $1\leq k\leq n-1$, 
then $f$ satisfies condition \eqref{1} at every point. Moreover, if inequality \eqref{lambda} 
is strict, then inequality \eqref{1} is also strict at every point.
\end{proposition}
\proof
Since 
$$
a(n,k,H,0)=\frac{n^2H^2}{n-k},
$$ 
the proof follows directly from the inequalities 
$S\leq n\lambda^2_n$ and $H\geq\lambda_1$.\qed
\vspace{1.5ex}

A large class of ellipsoids satisfies condition 
\eqref{lambda}. Consider, for instance, the 
ellipsoid in $\R^{n+1}$ defined by 
$$ 
\frac{x_1^2}{a_1^2}+\cdots+
\frac{x_{n+1}^2}{a_{n+1}^2} = 1,
$$
where $0<a_1\leq\dots\leq a_{n+1}$. A 
straightforward computation shows that the 
minimum and the maximum of the principal 
curvatures of the ellipsoid are $a_1/a_{n+1}^2$ 
and $a_{n+1}/a_1^2$, respectively. It follows 
that condition \eqref{lambda} holds provided 
$$
a_{n+1}\leq a_1\big({\frac{n}{n-k}}\big)^{1/6}. 
$$

By Hadamard’s classical theorem, every 
ovaloid in $\R^{n+1}$ is diffeomorphic to 
$\Sf^n$. Consequently, Propositions \ref{ex} 
and \ref{ell} provide geometrically distinct 
isometric immersions of manifolds 
diffeomorphic either to the sphere 
$\Sf^n$ or to the torus 
$\Sf^1 \times \Sf^{n-1}$ for $n \geq 4$, 
which also satisfy condition \eqref{1} with $c=0$ 
and $k=2$. Thus, we obtain numerous 
compact, geometrically distinct 
submanifolds that strictly satisfy 
\eqref{1} at every point. Moreover, 
this strict form is preserved under 
sufficiently small smooth deformations 
of any such example. Hence our 
results are optimal.

\begin{example}\emph{
Consider a rotational hypersurface in 
$\R^{n+1},n\geq4$, obtained by rotating 
a curve $c$ in the $x_1x_2$-plane 
around the $x_1$-axis. Suppose the 
curve is given as a graph of a positive 
function $u=u(x_1)$. After choosing 
the orientation appropriately, the 
rotational hypersurface has two 
principal curvatures given by
$$
\lambda=\frac{1}{u(1+(u')^2)^{1/2}},\;\,
\mu=-\frac{u''}{(1+(u')^2)^{3/2}},
$$
where $\lambda$ has multiplicity $n-1$. 
The hypersurface satisfies condition \eqref{1} 
with $k=2$ in strict form provided that 
$$
(n-1)\lambda^2+2(n-1)\lambda\mu-(n-3)\mu^2>0, 
$$
or equivalently, 
$$
-a_n(1+(u')^2)<uu''<b_n(1+(u')^2),
$$
where
$$
a_n=\frac{\sqrt{2(n-1)(n-2)}+n-1}{n-3},\;\,
b_n=\frac{\sqrt{2(n-1)(n-2)}-(n-1)}{n-3}.
$$
By connecting such hypersurfaces to spheres, 
as in Example 3.3 in \cite{MN}, we obtain 
compact hypersurfaces that satisfy condition 
\eqref{1} with $k=2$, whose first Betti number 
can be any nonnegative integer. The homology 
of these hypersurfaces is precisely as described 
in part (iii) of Theorem \ref{ck=2}.
}
\end{example}

\noindent Theodoros Vlachos\\
University of Ioannina \\
Department of Mathematics\\
45110 Ioannina -- Greece\\
e-mail: tvlachos@uoi.gr

\end{document}